\documentclass{icmart}

\contact[ifarah@yorku.ca]{}

\newtheorem{theorem}{Theorem}[section]

\newtheorem{lemma}[theorem]{Lemma}
\newtheorem{proposition}[theorem]{Proposition}

\newtheorem{question}[theorem]{Question}

\theoremstyle{definition}

\usepackage{amsmath, amssymb}
\usepackage{hyperref}
\usepackage{tikz}

\newcounter{my_enumerate_counter}
\newcommand{\pushcounter}{\setcounter{my_enumerate_counter}{\value{enumi}}}
\newcommand{\popcounter}{\setcounter{enumi}{\value{my_enumerate_counter}}}

\newcommand{\bfT}{\mathbf T} 
\newcommand{\bfP}{\mathbf P}

\DeclareMathOperator{\Irr}{Irr}
\DeclareMathOperator{\GL}{GL}

\DeclareMathOperator{\id}{id}
\newcommand{\cU}{\mathcal U}
\newcommand{\cO}{\mathcal O}

\newcommand{\bt}{\mathbf t}

\newcommand{\bbF}{{\mathbb F}}

\newcommand{\bbZ}{{\mathbb Z}}

\newcommand{\bft}{\mathbf t}

\newcommand{\bbK}{{\mathbb K}}

\newcommand{\bbN}{{\mathbb N}}
\newcommand{\bbC}{\mathbb C}

\newcommand{\bbR}{\mathbb R}

\newcommand{\cZ}{{\mathcal Z}}
\newcommand{\cA}{{\mathcal A}}

\newcommand{\fM}{\mathfrak M}

\newcommand{\rs}{\restriction}

\newcommand{\cC}{\mathcal C}

\newcommand{\cV}{\mathcal V}

\newcommand{\cB}{\mathcal B}

\newcommand{\cQ}{\mathcal Q}

\newcommand{\eps}{\varepsilon}

\DeclareMathOperator{\dist}{dist}

\newcommand{\cP}{\mathcal P} 
\DeclareMathOperator{\Th}{Th}
\DeclareMathOperator{\ThE}{Th_\exists}

\newcommand{\aet}{$\forall\exists$}

\newcommand{\snus}{separable, nuclear, unital and simple}

\DeclareMathOperator{\Sp}{sp}
\DeclareMathOperator{\Fin}{Fin}

\newcommand{\PNF}{\cP(\bbN)/\Fin}
 
\DeclareMathOperator{\Eq}{E}
\newcommand{\tfc}{2^{\mathfrak c}} 

\title[Logic and operator algebras]{Logic and operator algebras} 

\author[Ilijas Farah]
{Ilijas Farah\thanks{Partially supported by NSERC.}}

\begin{document}

\begin{abstract}
The most recent wave of applications of logic to operator algebras is a young and 
 rapidly developing field. This is 
 a snapshot of the current state of the art.  
   \end{abstract}



\maketitle


\section{Introduction}

The  connection between
logic and operator algebras in the past century was  
sparse albeit fruitful. 
Dramatic progress has brought set theory and operator algebras
closer together over the last decade. A number
of long-standing problems in the theory of C*-algebras were solved by
using set-theoretic methods, and solutions to some of them were even shown to be independent from ZFC. 
There is much to be said about these developments (as witnessed in 
three almost disjoint recent survey papers \cite{We:Set}, \cite{FaWo:Set}, \cite{Fa:Selected}), but that is not what this paper is about. 
New applications of logic to operator algebras are being found at such a pace that any survey 
is bound to become obsolete within a couple of years.
Instead of presenting an encyclopaedic survey, I shall proceed to describe the current developments 
(many of them from the unpublished joint work, \cite{FaHaRoTi}, \cite{Fa:MFO2013}) and outline some 
possible directions of research. The choice of the material reflects my interests and no attempt at  completeness 
has been made.  Several results proved by operator algebraists without using logic 
that have  logical content are also included.

`Logic' in the title refers to  model theory and (mostly descriptive) set theory, with a dash of recursion theory in a crucial place. 

I am indebted to Bradd Hart, Isaac Goldbring and Aaron Tikuisis
 for a number of remarks on the first draft of the present paper that considerably improved it in many ways.

\section{Operator algebras}
\label{S.OA} 

\label{S.GNS} 

Let $\cB(H)$ denote the Banach algebra of bounded linear operators on a complex Hilbert space $H$
equipp\-ed with the operation * of taking the adjoint. 
A \emph{C*-algebra} is a  Banach algebra with involution which is *-isomorphic to a subalgebra
of $\cB(H)$ for some $H$. Notably, all algebraic isomorphisms between C*-algebras are isometries.  
All C*-algebras considered here will be unital, unless otherwise specified. 
A \emph{von Neumann algebra} is a unital subalgebra of $\cB(H)$ which is closed in 
the weak operator topology. An algebra isomorphic to a von Neumann algebra is called a \emph{W*-algebra}. 
Standard terminology from operator theory is imported into operator algebras, 
and in particular positivity of self-adjoint operators plays an important role.

I only have something to say about those von Neumann algebras that have a trace. A  \emph{normalized trace} 
(on a von Neumann algebra or a unital C*-algebra) is a  unital positive functional  such that $\tau(ab)=\tau(ba)$ for all $a$ and $b$. 
We shall only consider unital algebras and normalized traces. 
A trace on  a von Neumann algebra is automatically continuous in the weak operator topology.  
A tracial infinite-dimensional von Neumann algebra with a trivial center is a \emph{II$_1$ factor}. 
The terminology comes from von Neumann's type classification, in which the unique I$_n$ factor
is $M_n(\bbC)$; we shall not consider other types of factors.

If $\tau$ is a trace on an operator algebra $A$ then the  $\ell_2$-norm  
$\|a\|_2=\tau(a^*a)^{1/2}$ turns $A$ into a pre-Hilbert space. 
The algebra $A$ is represented on this space by the left multiplication; 
this is the \emph{GNS representation} corresponding to $\tau$. 
If $A$ is a C*-algebra, then the weak closure of the image of $A$ is a tracial von Neumann algebra. 
If $A$ is simple and infinite-dimensional, this algebra is a II$_1$ factor. 
A GNS representation can be associated to an arbitrary positive unital functional (\emph{state}). 

The category of abelian C*-algebras is equivalent to the category of locally compact Hausdorff spaces 
and the category of abelian von Neumann algebras with a distinguished trace is equivalent to the category of 
measure algebras. Because of this, these two subjects are considered to be noncommutative (or quantized) 
topology and measure theory, respectively.

There is only one (obvious, spatial) way to define the tensor product of  von Neumann algebras. 
A C*-algebra $A$ is \emph{nuclear} if for every C*-algebra $B$ there is a unique C*-norm on the algebraic tensor product of $A$ and $B$. 
The importance of this notion is evident from a variety of its equivalent characterizations (see \cite{BrOz:C*}), one of them being 
Banach-algebraic amenability.  
Although by a result of Junge and Pisier (see \cite{BrOz:C*}) there is finite subset $F\Subset \cB(H)$ such that no nuclear C*-algebra 
includes $F$, these algebras are ubiquitous in a number of applications.  

For more on C*-algebras and von Neumann algebras see \cite{Black:Operator},  \cite{Jon:von}, and \cite{BrOz:C*}. 
\vskip.2in

\subsection{Intertwining} \label{S.Intertwining} A \emph{metric structure} is a complete metric space $(A,d)$
equip\-ped with functions $f\colon A^n\to A$ and predicates $p\colon A^n\to \bbR$, all of which are assumed to be 
uniformly continuous on $d$-bounded sets.  
Consider  
two separable complete metric structures $A$ and $B$. 
Assume  we have partial isometric homomorphisms $\Phi_n\colon F_n\to G_n$, $\Psi_n\colon G_n\to F_{n+1}$  for $n\in \bbN$ 
 such that $F_n\subseteq F_{n+1}\subseteq A$ and 
$G_n\subseteq G_{n+1}\subseteq B$ for $n\in \bbN$ and  $\bigcup_n F_n$ and $\bigcup_n G_n$ are dense in $A$ and $B$ respectively. 
Furthermore assume that  in the following diagram

\begin{tikzpicture}
  \matrix[row sep=1cm,column sep=1cm] {
\node  (A1) {$F_1$};
& \node (A2) {$F_2$};
& \node (A3) {$F_3$};
 & \node (A4) {$F_4$};
& \node (Adots) {$\dots$};
& \node (A) {$A$};
\\
\node  (B1) {$G_1$};
& \node (B2) {$G_2$};
& \node (B3) {$G_3$};
 & \node (B4) {$G_4$};
& \node (Bdots) {$\dots$};
& \node (B) {$B$};
\\
};
\draw (A1) edge[->] node [left] {$\Phi_1$} (B1) ;
\draw (A2) edge[->] node  [left] {$\Phi_2$} (B2) ;
\draw (A3) edge[->] node [left] {$\Phi_3$} (B3);
\draw (A4) edge[->] node [left] {$\Phi_4$} (B4);
\draw (B1) edge[->] node [left] {$\Psi_1$} (A2);
\draw (B2) edge[->] node [left] {$\Psi_2$} (A3);
\draw (B3) edge[->] node [left] {$\Psi_3$}  (A4);
\draw (B4) edge[->] (Adots);
\draw (A1) edge[->] (A2);
\draw (A2) edge[->] (A3);
\draw (A3) edge[->] (A4);
\draw (A4) edge[->] (Adots);
\draw (B1) edge[->] (B2);
\draw (B2) edge[->] (B3);
\draw (B3) edge[->] (B4);
\draw (B4) edge[->] (Bdots);

\end{tikzpicture}

\noindent the $n$-th triangle commutes up to $2^{-n}$. 
Then $\Phi\colon \bigcup_n F_n\to B$ defined by $\Phi(a)=\lim_n \Phi_n(a)$ and $\Psi\colon \bigcup_n G_n\to A$
defined by $\Psi(b)=\lim_n \Psi_n(b)$ are well-defined isometric homomorphisms.  Their continuous extensions to $A$ and $B$ 
are respectively  an isomorphism from $A$ onto $B$ and vice versa. 

Variations of this method for constructing isomorphisms between C*-algebras comprise  \emph{Elliott's intertwining argument}.
In \emph{Elliott's program} for classification of \snus{} C*-algebras maps $\Phi_n$ and $\Psi_n$ are obtained by lifting morphism
between the K-theoretic invariants (so-called \emph{Elliott invariants}) of $A$ and $B$. 
The first result along these lines was the Elliott--Bratteli classification of separable AF algebras (i.e., direct limits 
of finite-dimensional C*-algebras) by the ordered $K_0$.  
Remarkably, for $A$ and $B$ belonging to a rather large class of  nuclear  C*-algebras this method shows that 
any morphism between Elliott invariants lifts to a morphism between the algebras. 
Elliott conjectured that the \snus{} algebras are classified by K-theoretic invariant known as the Elliott invariant. 
This bold conjecture was partially confirmed in many instances.  
See \cite{Ror:Classification} for more on the early history of this 
fascinating subject. 

Examples of \snus{} C*-algebras that limit the extent of Elliott's classification program 
 were given in \cite{Ror} and \cite{To:On}. 
Algebras defined in  \cite{To:On}
have a remarkable additional property. Not only do the nonisomorphic algebras $A$ and $B$ 
have the same Elliott invariant, but in addition they cannot be distinguished 
by any homotopy-invariant continuous functor. We shall return to these examples in \S\ref{S.Crazy}. 
The revised Elliott program is still one of the core subjects in 
the study  of C*-algebras (see \cite{EllTo:Regularity}). 

\vskip.2in
\subsection{Strongly self-absorbing (s.s.a.) algebras} \label{S.s.s.a.}
An infinite-di\-men\-sio\-nal  C*-algebra is \emph{UHF} (uniformly hyperfinite) if it is an infinite tensor product of full matrix algebras $M_n(\bbC)$. 
If $A$ is UHF, then every two unital copies of~$M_n(\bbC)$ in it are unitarily conjugate and therefore 
every endomormphism  of $A$ is a point-norm limit of inner automorphisms. 
The \emph{generalized natural number} of $A$ has as its `divisors' 
all $n$ such that $M_n(\bbC)$ embeds unitally into $A$. Glimm proved that this is a complete isomorphism 
invariant for the separable UHF algebras. 
 
 If $A$ is UHF then it has a unique trace $\tau$. The tracial von Neumann algebra corresponding to the $\tau$-GNS representation 
 of $A$ (\S\ref{S.GNS}) is the \emph{hyperfinite II$_1$ factor, $R$}, and it does not depend on the choice of $A$. 
 It is the only injective II$_1$ factor and it has played a key role in the classification of injective factors (\cite{Connes:Class}). 
 
Two *-homomorphisms $\Phi$ and $\Psi$ from $A$ into $B$ are \emph{approximately unitarily equivalent} 
if there is a net of inner automorphisms $\alpha_\lambda$, for $\lambda\in \Lambda$, of $B$ such that 
\[
\textstyle \lim_\lambda \alpha_\lambda\circ \Phi(a) =\Psi(a)
\]
for all $a\in A$ (convergence is taken in the operator norm for C*-algebras and in the 
$\ell_2$-norm for tracial von Neumann algebras).   
If $A\otimes B\cong A$   we say that $A$ is \emph{$B$-absorbing} and if $A\otimes A\cong A$ then we say that 
$A$ is  \emph{self-absorbing}.  
Here and in what follows, we will often be providing two definitions at once, one
for von Neumann algebras and another for C*-algebras. The difference
comes in the interpretation of $\otimes$, either as the von Neumann
(spatial) tensor product $\bar\otimes$ or as the C*-algebra minimal (spatial) tensor product $\otimes$.
  \emph{McDuff} factors are the $R$-absorbing II$_1$ factors. 
A separable C*-algebra $D$ is \emph{strongly self-absorbing} (s.s.a.) (\cite{ToWi:Strongly}) if 
there is an isomorphism $\Phi\colon D\to D\otimes D$ and 
map $\id\otimes 1_D:D\to D\otimes D$ is approximately unitarily equivalent to $\Phi$. 
The definition of strongly self-absorbing is modified to II$_1$ factors following the 
convention stated above, by replacing $\|\cdot\|$ with $\|\cdot\|_2$ 
and $\otimes$ with $\bar\otimes$. 

The hyperfinite factor  $R$ is the only s.s.a. tracial von Neumann algebra with separable predual 
(Stefaan Vaes pointed out that this was essentially proved in  \cite[Theorem 5.1(3)]{Connes:Class}). 
A UHF algebra $A$ is s.s.a. if and only if it is self-absorbing. 
However, the latter notion is in general much stronger. For any 
unital C*-algebra $A$ the infinite tensor product $\bigotimes_{\bbN} A$
is self-absorbing but not necessarily s.s.a. 
Every  s.s.a. C*-algebra $D$ is  simple, nuclear and unital~(\cite{effros1978c}).

Three s.s.a. algebras are particularly important. 
The \emph{Jiang--Su} algebra $\cZ$ is an infinite-dimensional C*-algebra which is indistinguishable from $\bbC$ 
by its Elliott invariant. Conjecturally, $\cZ$-absorbing infinite-dimensional \snus{} algebras are classifiable by their Elliott invariant. 
The Cuntz algebra~$\cO_2$ is the universal algebra generated by two partial isometries with complementary ranges. 
The Cuntz algebra $\cO_\infty$ is the universal unital C*-algebra generated by partial isometries $v_n$, for $n\in \bbN$, 
with orthogonal ranges.  
The first step in the Kirchberg--Phillips  classification of purely infinite \snus{} algebras was Kirchberg's result that 
every such algebra is $\cO_\infty$-absorbing  and that $\cO_2$ is $A$-absorbing for every \snus{} algebra
 (see \cite{Ror:Classification}).

\section{Abstract classification}
\label{S.Abstract}

A \emph{Polish space} is a separable, completely metrizable topological space. 
A subset of a Polish space is \emph{analytic} if it is a continuous image of some Polish space. 
Essentially all classical classification problems in mathematics (outside of subjects with a strong 
set-theoretic flavour)  can be modelled by an analytic equivalence 
relation on a Polish space. Moreover, the space of classifying invariants is also of this form, and
computation of the invariant is usually given by a Borel measurable map. 
This is indeed the case with C*-algebras and the Elliott invariant~(\cite{FaToTo:Descriptive}). 

If $E$ and $F$ are equivalence relations on Polish spaces, $E$ is \emph{Borel-reducible} to $F$, $E\leq_B F$, 
if there exists a Borel-measurable $f\colon X\to Y$ such that $x\Eq y$ if and only if $f(x)\Eq f(y)$. 
One can interpret this as stating that the classification problem for $E$ is not more difficult than the classification problem for $F$. 
Following Mackey, an equivalence relation $E$  Borel-reducible to the equality relation on some Polish space is said to be \emph{smooth}.   
By the \emph{Glimm--Effros} dichotomy the class of non-smooth Borel-equivalence relations has an initial object
(\cite{HaKeLo:GE}), denoted $E_0$. It is the tail equality relation on $\{0,1\}^{\bbN}$. 
 While the Glimm--Effros dichotomy was proved by using sophisticated tools from effective descriptive set theory, the combinatorial core of the proof
can be traced back to work of Glimm and Effros on representations of locally compact groups and separable C*-algebras. 
See \cite{Ke:Classical}, \cite{Gao:Invariant} for more on (invariant) descriptive set theory.

When is an equivalence relation classifiable? 
Many non-smooth equivalence relations are considered to be satisfactorily classified. 
An example from the operator algebras is  the Elliott--Bratteli classification of separable AF algebras
 by countable abelian ordered groups. A rather generous notion is being `classifiable by countable structures.' 
 Hjorth's theory of turbulence (\cite{Hj:Book})
provides a powerful tool for proving that an orbit equivalence relation is not classifiable by countable structures.

Sasyk and T\"ornquist    
have 
proved that every class of injective factors that was not already satisfactorily classified is not classifiable by countable structures (\cite{sato09b}, \cite{sasyk2010turbulence}). 
  By combining results of \cite{FaToTo:Turbulence}, \cite{Mell:Computing}, \cite{GaKe}, \cite{EFPRTT}  and \cite{Sab:Completeness},
 one proves that the following isomorphism relations are Borel-equireducible.\\ 
(a)  Isomorphism relation of separable C*-algebras. \\
(b)  Isomorphism relation of Elliott--classifiable \snus{}  algebras.\\ 
(c) Isometry relation of separable Banach spaces.\\ 
(d)  Affine homeomorphism relation of metrizable Choquet simplices. \\
(e)  Isometry relation of Polish spaces. 

Each of these equivalence relations  (as well as the isometry of a class of separable metric 
structures of any given signature) is   Borel-reducible to an orbit equivalence 
relation of  a Polish group action (\cite{EFPRTT}).

Being Borel-reducible to an orbit equivalence relation is, arguably, the  most generous definition
 of being concretely classifiable. Conjecturally, $E_1$, the tail-equivalence relation on $[0,1]^{\bbN}$, 
 is an initial object  among Borel equivalence relations not Borel-reducible to an orbit equivalence relation (\cite{KeLou:Structure}). 
Notably, the isomorphism of separable Banach spaces
 is the $\leq_B$-terminal object among analytic equivalence relations (\cite{FeLouRo}). 

The answer to the  question `When is an equivalence relation classifiable' is frequently of somewhat sociological nature. 
 It is notable that the isomorphism relation of abelian unital C*-algebras (generally considered intractable) 
is Borel-reducible to the isomorphism relation of Elliott-classifiable AI  algebras (for which there is a satisfactory classification relation).  
Also, as pointed out by David Fremlin, most analysts find that  normal operators are satisfactorily classified up to conjugacy by 
the spectral theorem, although they are not classifiable by countable structures. 

Nevertheless, the theory of Borel-reducibility is a great example of a situation in which logic provides concrete obstructions to 
sweeping conjectures. 
For example, 
 the classification of countable abelian torsion free groups of rank $n+1$ is strictly more 
complicated than the classification of countable abelian torsion free groups of rank $n$ for every $n$ (\cite{Tho:On}). 
(Notably, the proof of this result uses \emph{Popa superrigidity} of  II$_1$ factors, \cite{Po:Strong}.)
This theory was recently successfully applied to (non)classification of automorphisms of group actions 
on operator algebras (\cite{KeLiPi:Turbulence},  automorphisms of C*-algebras (\cite{Lup:Unitary}, \cite{KeLuPh:Borel})
and subfactors (\cite{brothier2013families}).

   A partial   Borel-reducibility diagram of classification problems in operator algebras is 
given below. For an explanation of terminology see \cite[\S 9]{Fa:Selected}.  
I am indebted to Marcin Sabok for pointing out that the isomorphism of countable structures of any signature 
is Borel-reducible to the isomorphism relation of separable AF algebras
(\cite{camerlo2001completeness}).

{\tiny
\usetikzlibrary{shapes}

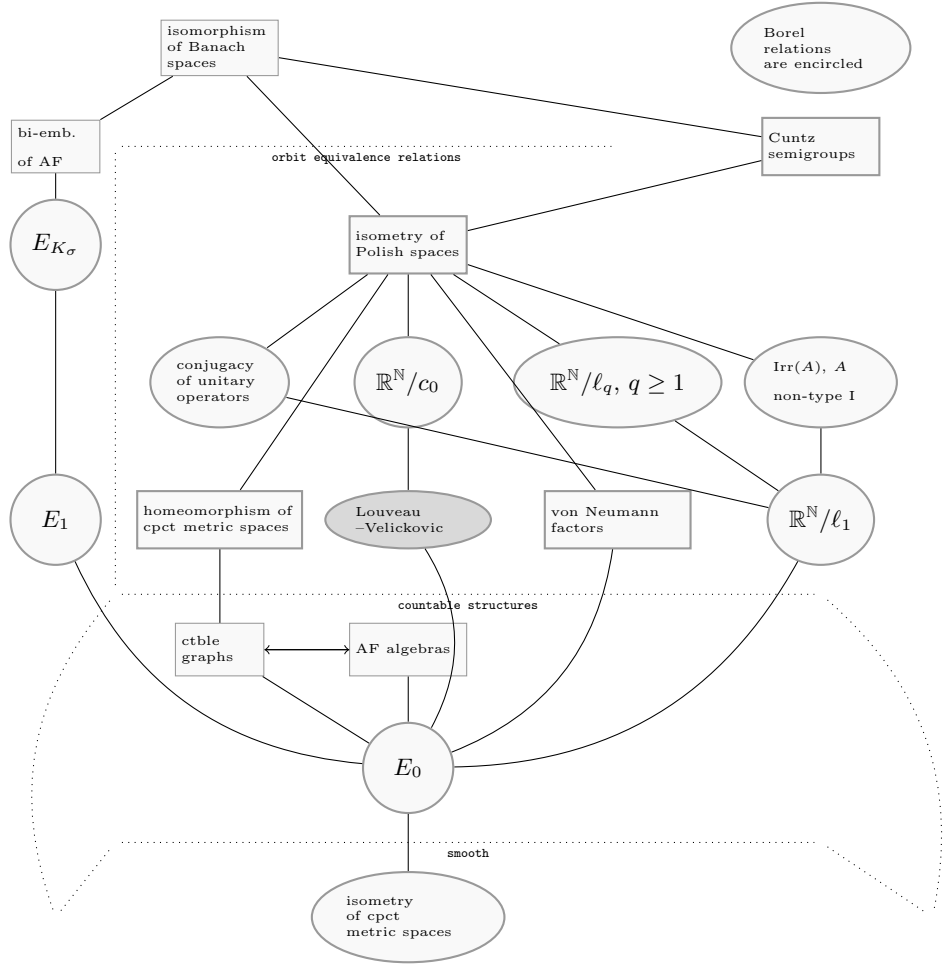
\tikzstyle{smooth}=[ellipse,
										font={\small},
                                                thick,
                                                text centered,
                                                minimum size=1.2cm,
                                                draw=gray!80,
                                                fill=gray!05]

\tikzstyle{dark-side}=[rectangle,
										font={\small},
                                                thin,
                                                minimum size=.7cm,
                                                draw=gray!80,
                                                fill=gray!05]

\tikzstyle{unknown}=[rectangle,
										font={\small},
                                                thick,
                                                minimum size=.3in,
                                                draw=gray!80,
                                                fill=gray!05]

\tikzstyle{orbit}=[rectangle,
										font={\small},
                                                thick,
                                                minimum size=.3in,
                                                draw=gray!80,
                                                fill=gray!05]

\tikzstyle{big-smooth}=[ellipse,
										font={\small},
                                                thick,
                                                minimum size=.3in,
                                                draw=gray!80,
                                                fill=gray!30]

\tikzstyle{ctble}=[rectangle,
										font={\small},
                                                thick,
                                                minimum size=1.2cm,
                                                draw=gray!80,
                                                fill=gray!05]

{\tiny \noindent\begin{tikzpicture}
  \matrix[row sep=0.3cm,column sep=0.1cm] {
&& &
\node  (Banach)[dark-side] {\parbox{.55in}{\tiny isomorphism\\  of Banach\\ spaces}}; &
&&  \node (note)[smooth]{\parbox{.6in}{\tiny Borel\\ relations\\ are encircled}};
\\
 \node (biAF)[dark-side] {\parbox{.4in}{{\tiny bi-emb.\\ of AF}}};
& \node (orbit-l) {};
& &  & &
   \node (orbit-r) {};  &
\node  (Cuntz)[unknown] {\parbox{.55in}{\tiny Cuntz\\ semigroups}};
\\
\node (EKsigma)[smooth] {$E_{K_\sigma}$};
& & &
&
\node (iso-Polish)[orbit] {\parbox{.55in}{\tiny isometry of\\ Polish spaces}};
& &
\\
&& &
&
&
\\
&
&
& \node(conjugacy)[smooth] {\parbox{.45in}{\tiny conjugacy\\ of unitary\\ operators}};
& \node (czero) [smooth] {$\bbR^\bbN/c_0$};
&  \node(ellq)[smooth]{$\bbR^\bbN/\ell_q$, $q\geq 1$};
&\node(Irr)[smooth]{\parbox{.5in}{{\tiny $\Irr(A)$, $A$\\ non-type I}}};
\\
& &  & 
&
&
& 
&
\\
\node (Eone)[smooth]{$E_1$};
&&
& \node  (cpct)[orbit]{\parbox{.8in}{\tiny homeomorphism  of\\ cpct metric spaces}};
& \node (LV) [big-smooth] {\parbox{.55in}{\tiny Louveau\\ --Velickovic}};
&
\node (vNA)[orbit] {\parbox{.7in}{\tiny von Neumann\\ factors}};
& \node (ell1) [smooth] {$\bbR^\bbN/\ell_1$};
&
\\
& \node (ctble-l) {};
& & &  & & \node  (ctble-r) {};
\\
& &  &  \node (graph)[dark-side]{\parbox{.4in}{\tiny ctble\\ graphs}} ;
  &   
 \node (AF) [dark-side] {\parbox{.55in}{\tiny AF algebras}};
& 
\\
&  & & 
&
\\
& & & &
\node (Ezero) [smooth] {$E_0$ };
&
&
\\
&  \node  (smooth-l) {};
& & &  & &
\node (smooth-r) {};
\\
\node (bottom-l) {};
& &
&& 
\node  (cpct-isometry)[smooth]{\parbox{.65in}{ \tiny isometry\\ of cpct\\ metric spaces}};
& 
&
&&
\node (bottom-r) {};
\\
};
\path (Banach)edge  (biAF);
\path (conjugacy) edge (iso-Polish);
\path [double] (AF) edge [<->]  (graph);
\path (Ezero) edge [bend right] (vNA);
\path (Ezero) edge (AF);
\path (Ezero) edge (graph);
\path (ell1) edge  (conjugacy);
\path (Irr) edge (iso-Polish);
\path (AF) edge (graph);
\path (ellq) edge  (iso-Polish);
\path (ell1) edge (ellq);
\path (ell1) edge (Irr);
\path (Ezero) edge [bend right]   (ell1);
\path (czero) edge  (iso-Polish);
\path (Ezero) edge [bend right]   (LV);
\path (czero) edge    (LV);
\path (graph) edge (cpct);
\path (Eone) edge  [bend right] (Ezero);
\path (Ezero) edge (cpct-isometry);
\path (Eone) edge (EKsigma);
\path (EKsigma) edge (biAF);
\path (iso-Polish) edge  (Cuntz);
\path(vNA) edge (iso-Polish);
\path(iso-Polish) edge (Banach);
\path (Cuntz) edge (Banach);
\path  (cpct) edge (iso-Polish);
\path[dotted] (orbit-l) edge node[below]{\tt orbit equivalence relations} (orbit-r);
\path[dotted](orbit-l) edge (ctble-l);
\path[dotted](ctble-l) edge [bend right] (bottom-l);
\path[dotted](ctble-r) edge [bend left] (bottom-r);
\path[dotted] (ctble-l) edge node[below]{\tt countable structures} (ctble-r);
\path[dotted] (smooth-l) edge node[below]{\tt smooth} (smooth-r);
\path[dotted](smooth-r) edge  (bottom-r);
\path[dotted] (bottom-l) edge (smooth-l);
\end{tikzpicture}
}}

Borel-reduction of equivalence relations as defined above does not take into the account 
the functorial  nature of the classification of C*-algebras. Some preliminary results on Borel functorial 
classification were obtained by Lupini.

\section{Model-theoretic methods} 

Until recently there was not much interaction between model 
theory and  operator algebras (although model theory was fruitfully applied to the geometry 
of Banach spaces, see \cite{hensonultraproducts}). 
 Recent  emergence of the logic of metric structures~(\cite{BYBHU}), originally introduced only for bounded metric structures, 
created new opportunities  for such interactions. 
It was modified to allow  operator algebras 
in~\cite{FaHaSh:Model2}.

\subsection{Logic of metric structures} 
Model theory can roughly be described as  the study of axiomatizable classes of structures and sets definable in them. 
Axiomatizable properties can be  expressed in syntactic terms,  but they are also  characterized by  
preservation under ultraproduts and ultraroots (see \S\ref{S.Ultra}).  
A category~$\cC$ is \emph{axiomatizable} if there exists  a first-order theory $\bfT$  such that 
the  category $\fM(\cC)$ of all models of $\bfT$ is equivalent to the original category.

Classical model theory deals with discrete structures, and 
its variant suitable for metric structures as defined in \S\ref{S.Intertwining} was introduced in \cite{BYBHU}. 
In this logic interpretations of formulas are real-valued, propositional connectives are real-valued functions, and 
quantifiers  are $\sup_x$ and $\inf_x$. Each function  and predicate symbol is equipped with a modulus of uniform 
continuity. This modulus is a part of the language. 
If the diameter of the metric structures is fixed, then every formula has its own modulus of uniform continuity, 
respected in all relevant metric structures.   
Formulas form a real vector space  
equipped with a seminorm, $\|\phi(\bar x)\|=|\sup \phi(\bar a)^A|$ where
the supremum is taken over all metric structures $A$ of the given language and all tuples $\bar a$ in $A$ of the appropriate type. 
Formulas are usually required to have range in $[0,\infty)$ (or $[0,1]$ in the bounded case) but allowing negative 
values results in equivalent logic; see also \cite{BY:On}. 
The \emph{theory} of a model is the kernel of the functional $\phi\mapsto \phi^A$, where $\phi$ ranges over all sentences (i.e., formulas with no free variables)  of the language. This kernel uniquely defines the functional, which can alternatively be identified with the theory.  
The weak*-topology on this space is also known as the \emph{logic topology}. If the language is countable 
then the  space of formulas is separable and the spaces of theories and types (see  \S\ref{S.Omitting}) are equipped
with compact metric topology. 

Two metric structures are \emph{elementarily equivalent} if their theories coincide.   
A formula is \emph{existential} if it is of the form $\inf_{\bar x}\phi(\bar x)$ for some quantifier-free formula~$\phi(\bar x)$. 
The \emph{existential theory} of $A$ is $\ThE(A)=\{\psi\in \Th(A): \psi$ is existential$\}$.

There are several equivalent ways to adapt the  logic of metric structures
 to operator algebras  and to unbounded metric structures in general (\cite{FaHaSh:Model2},
 \cite{BY:Continuous}). 
Axiomatizability is defined via equivalence of categories as above, but
 model $M(A)$ associated with  $A$ has more (albeit artificial) structure.
 It is equipped with \emph{domains of quantification}, 
 bounded subsets of $A$ on which all functions and predicates are uniformly continuous (with a fixed modulus of uniform continuity) and 
 over which quantification is allowed. 
  It is the existence of category $\fM(\cC)$, and not its particular  choice, that matters.

  In the simplest  version of $M(A)$ quantification is allowed only over the (operator norm) $n$-balls  of the algebra. 
  The notion of \emph{sorts} over which one can quantify corresponds to those functors from the model 
  category into metric spaces with uniformly continuous functions that commute with ultraproducts~(see \cite[\ 2]{FaHaSh:Model2}). 
For example,      $M(A)$ can be taken  to consist of all matrix algebras $M_n(A)$  
 for $n\in \bbN$, as well as completely positive, contractive maps between them and finite-dimensional algebras. 
 This is important because nuclearity is equivalently characterized as the CPAP, the \emph{completely positive approximation property}
  (see \cite{BrOz:C*} and  \cite{Bro:Symbiosis}). 
 
C*-algebras are axiomatized as  Banach algebras with an involution that satisfy the \emph{C*-equality}, $\|aa^*\|=\|a\|^2$, 
by  the Gelfand--Naimark and Segal (GNS mentioned earlier) theorem. 
Abelian C*-algebras are  obviously axiomatized by $\sup_{x,y} \|xy-yx\|$ and 
non-abelian C*-algebras are slightly less obviously  axiomatized by $\inf_{\|x\|\leq 1} |1-\|x\||+\|x^2\|$  
(a C*-algebra is nonabelian if and only if it contains a nilpotent element).

The proof that the  tracial von Neumann algebras with a distinguished trace are also axiomatizable 
(\cite{FaHaSh:Model2}, first proved in \cite{BYHJR}) goes deeper and 
uses Kaplansky's Density Theorem. Again, quantification is allowed over the (operator norm) unit ball 
and the metric is the $\ell_2$ metric $\|a\|_2=\tau(a^*a)^{1/2}$. 
The operator norm is not continuous with respect to the $\ell_2$ metric and it therefore cannot be added to 
II$_1$ factors as a predicate. 

There are elementarily equivalent but nonisomorphic separable unital AF algebras. 
This is proved by using descriptive set theory. 
The association $A\mapsto \Th(A)$ is Borel, and  hence the relation of elementary equivalence is 
smooth (\S\ref{S.Abstract}). 
The category of AF algebras is equivalent to the category of their ordered $K_0$ groups. 
By the Borel version of this result and the fact that the isomorphism of dimension groups is not smooth 
 the conclusion follows. 

The following  proposition is taken from  (\cite{Fa:MFO2013}). 

\begin{proposition}\label{P.4.1}  
\begin{enumerate}
\item For every \snus{} C*-algebra there exists an elementarily 
equivalent, separable,  non-nuclear, C*-algebra. 
\item The reduced group C*-algebra  of the free group with infinitely many 
generators C$^*_r(F_\infty)$ is not elementarily 
equivalent to a nuclear C*-algebra. 
\end{enumerate}
\end{proposition}

Instead of providing a genuine  obstruction, this proposition precipitated some of the most interesting progress in the field.

Here is a simple but amusing observation.  
The \emph{Kadison--Kastler} distance between  subalgebras of $\cB(H)$ is the Hausdorff (norm) distance between their unit balls. 
 For every sentence $\phi$ the map $A\mapsto \phi^A$ is continuous with respect to this metric. 
 Therefore the negation of an axiomatizable property is stable under small perturbations of an algebra
 (see \cite{CSSWW:Perturbations} and references thereof for more on perturbations of C*-algebras).

\subsection{Elementary submodels} 
\label{S.ES}
If $A$ is a submodel of $B$, it is said to be an \emph{elementary submodel} if $\phi^A=\phi^B\rs A^n$ for every $n$ and every $n$-ary formula
$\phi$.   
The Downwards L\"owenheim--Skolem theorem implies that every model has a separable elementary submodel. 
Its C*-algebraic variant is known as `Blackadar's method' and is  used to provide separable
examples from known nonseparable examples (see \cite[II.8.5]{Black:Operator} and \cite{Oza:Dixmier}). 

\begin{proposition} \label{P.Elementary} 
Assume $A$ is a C*-algebra and $B$ is its elementary submodel. 
Then $B$ is a C*-algebra with the following properties. 
\begin{enumerate}
\item Every trace of $B$ extends to a trace of $A$. 
\item Every ideal of $B$ is of the form $I\cap B$ for some ideal $I$ of $A$. 
\item Every character of $B$ extends to a character of $A$. 
\item If $A$ is nuclear so is $B$. 
\end{enumerate}
\end{proposition} 

In particular, $B$ is monotracial and/or simple if and only if $A$ has these properties. 
It should be noted that neither of these properties is axiomatizable, because neither 
 of them is preserved under taking ultrapowers (see \cite{robert2013nuclear} for the 
nonaxiomatizability of having a unique trace and \S\ref{S.Ultra} for the ultrapowers).

A drastic example of a property that does not persist to elementary submodels is given in Theorem~\ref{T.amenable}. 

\subsection{Intertwining again} \label{S.Crazy}
We return to Elliott's intertwining argument~(\S\ref{S.Intertwining}): 

\begin{tikzpicture}
  \matrix[row sep=1cm,column sep=1cm] {
\node  (A1) {$A_1$};
& \node (A2) {$A_2$};
& \node (A3) {$A_3$};
 & \node (A4) {$A_4$};
& \node (Adots) {$\dots$};
& \node (A) {$A=\lim_n A_n$};
\\
\node  (B1) {$B_1$};
& \node (B2) {$B_2$};
& \node (B3) {$B_3$};
 & \node (B4) {$B_4$};
& \node (Bdots) {$\dots$};
& \node (B) {$B=\lim_n B_n$};
\\
};
\draw (A1) edge[->] node [left] {$\Phi_1$} (B1) ;
\draw (A2) edge[->] node  [left] {$\Phi_2$} (B2) ;
\draw (A3) edge[->] node [left] {$\Phi_3$} (B3);
\draw (A4) edge[->] node [left] {$\Phi_4$} (B4);
\draw (B1) edge[->] node [left] {$\Psi_1$} (A2);
\draw (B2) edge[->] node [left] {$\Psi_2$} (A3);
\draw (B3) edge[->] node [left] {$\Psi_3$}  (A4);
\draw (B4) edge[->] (Adots);
\draw (A1) edge[->] (A2);
\draw (A2) edge[->] (A3);
\draw (A3) edge[->] (A4);
\draw (A4) edge[->] (Adots);
\draw (B1) edge[->] (B2);
\draw (B2) edge[->] (B3);
\draw (B3) edge[->] (B4);
\draw (B4) edge[->] (Bdots);

\end{tikzpicture}

\noindent If the maps $\Phi_n$ are expected to converge to an isomorphism, it is necessary that 
they approximate elementary maps. For a formula $\phi(\bar x)$ and a tuple $\bar a$ in the domain of $\Phi_n$ one
must have $\phi(\bar a)^A=\lim_n \phi(\Phi_n(\bar a))^B$. 
Even more elementarily, the algebras $A$ and $B$ ought to be elementarily equivalent (no pun intended).  
Every known counterexample to Elliott's program involves \snus{} algebras
with the same Elliott invariant, but different theories.
For example, the radius of comparison was used in \cite{To:Comparison} to distinguish between continuum many nonisomorphis \snus{}
algebras with the same Elliott invariant, and it can be read off from the theory of an algebra  (\cite{Fa:MFO2013}).

This motivates an outrageous conjecture, that the following question has a positive answer. 

\begin{question} \label{Q.Crazy} Assume that  \snus{} algebras $A$ and $B$ have the  same Elliott invariant
and are elementarily equivalent. Are $A$ and $B$ necessarily isomorphic? 
\end{question} 

Since being $\cZ$-stable is axiomatizable (see \S\ref{S.Relative}), the revised Elliott conjecture that 
all $\cZ$-stable \snus{} algebras are classified by their Elliott invariant is a special case of a positive answer to Question~\ref{Q.Crazy}. 
All known nuclear C*-algebras belong to the so-called \emph{bootstrap class}, obtained by closing the class of type I algebras under 
operations known to preserve nuclearity (see~\cite{Black:Operator}). 
An (expected) negative answer to  Question~\ref{Q.Crazy}  would require new examples of \snus{} algebras. 
Can model-theoretic methods provide such examples?

\subsection{Omitting types} 
\label{S.Omitting} 

Let $\bbF_n$ be the set of formulas whose free variables are included in $\{x_1,\dots, x_n\}$. 
An \emph{$n$-type} is a subset $\bft$ of $\bbF_n$  such that for every finite $\bft_0\Subset \bft$
and for every $\eps>0$ there are a C*-algebra $A$ and $n$-tuple $\bar a$ in the unit ball of $A$ such that 
$|\phi(\bar a)|<\eps$ for all $\phi\in \bft_0$. By applying functional calculus one sees that this definition is equivalent 
to the apparently more general standard definition (\cite{BYBHU}) in which types consist of arbitrary closed conditions. 
An $n$-type $\bft$ is \emph{realized} by $n$-tuple $\bar a$  in the unit ball of $A$ 
if $\phi(\bar a)^A=0$ for all $\phi$ in $\bft$. It is \emph{omitted} in $A$ if it is not realized
by any $n$-tuple in $A_1$. 
\L o\'s's theorem implies that every type is realized in an ultraproduct; we shall return to this in \S\ref{S.CtbleSat}
but presently we are concerned with omitting types.

The omitting types theorem of classical (`discrete') model theory 
 (\cite{Mark:Model})  provides a simple condition for omitting a type  in a model of a given theory. 
A predicate $p$ is \emph{definable} if for every $\eps>0$ there exists a formula $\phi(\bar x)$ 
which up to $\eps$ approximates the value of $p$. A type is definable if the distance function to its realization
in a saturated model (\S\ref{S.CtbleSat}, \S\ref{S.Automorphisms}) is definable. 
By the omitting types theorem of \cite{BYBHU} a type is omissible if and only if it is not definable, with the additional 
stipulation that it be \emph{complete} (i.e., maximal under the inclusion). 
While a definable type is never omissible even if it is incomplete, Ben Yaacov has isolated
types that are neither definable nor omissible. His example was simplified by T. Bice.

\begin{theorem}[\cite{FaMa:Omitting}]\label{T.Types} 
\begin{enumerate}
\item There is a theory $\bfT$ in a separable language such that the set of types omissible in some model of $\bfT $
is a complete $\Sigma^1_2$ set. 
\item There are a complete theory $\bfT$ in a separable language and a countable set $\bfP$ of types 
such that for every finite $\bfP_0\Subset \bfP$ there exists a model 
$M$ of $\bfT$ that omits all types in $\bfP_0$, but no model of $\bfT$ 
omits all types in $\bfP$. 
\end{enumerate}
\end{theorem}

Therefore the question of whether a type is omissible in a model of a given metric theory
is by (1) far from being Borel or even analytic and therefore  intractable, 
and by (2) separately omissible types over a complete theory 
are not necessarily jointly omissible. Both results stand in stark contrast to the situation in 
classical model theory.

The idea that the omitting types theorem can be used  in the study of C*-algebras 
emerged independently in \cite{Scow:Some} and \cite{Mitacs2012}. 
A sequence $\bt_n$, for $n\in \bbN$, of $m$-types is \emph{uniform} if there are 
formulas $\phi_j(\bar x)$ for $j\in \bbN$ with the same modulus of uniform continuity 
such that $\bt_n=\{\phi_j(\bar x)\geq 2^{-n}: j\in \bbN\}$
for every $n$. 
In this situation, the interpretation of the infinitary formula $\phi(\bar x)=\inf_j \phi_j(\bar x)$ is uniformly continuous in 
every model (with a fixed modulus of uniform 
continuity) and moreover $\sup_{\bar x}\phi(\bar x)^A=0$ if and only if $A$ omits all $\bt_n$. 

Nuclearity, simplicity, as well as many other important non-axiomatizable properties of C*-algebras
(including nuclear dimension or decomposition rank $\leq n$; see \cite{WinterZacharias:NucDim})
are characterized by omitting a uniform sequence of types. 
 The classical theory of omitting types 
applies to such types unchanged:  a uniform sequence of types is omissible in a model of a complete theory $\bfT$ if and only if none of the types is isolated (\cite{FaMa:Omitting}).  
As an extra, this characterization shows that one can find a separable elementary submodel of a nonnuclear algebra that is itself nonnuclear by assuring that it includes a tuple that realizes the relevant type

\subsection{Strongly self-absorbing  algebras II} \label{S.s.s.a..model} 
These algebras have remarkable model-theoretic properties. 
 Every s.s.a. algebra $D$ is a  \emph{prime model} of its theory (it elementarily embeds into every other 
 model of its theory) and 
 every unital morphism of  $D$ into another model of its theory 
 is elementary~(\S\ref{S.ES}).

\begin{proposition} \label{P.ssa} 
If $D$ and $E$ are s.s.a. algebras then the following are equivalent. 
\begin{enumerate}
\item $E$ is \emph{$D$-absorbing}: $E\otimes D\cong E$. 
\item $D$ is isomorphic to a subalgebra of $E$. 
\item $\Th_\exists(D)\subseteq \Th_\exists(E)$. 
\end{enumerate}
\end{proposition} 

The implications from (1) to (2) and from (2) to (3) are always true, but both converses fail in general.  
S.s.a. algebras are as rare as they are important and  
the following diagram represents all known s.s.a. algebras, given in the order defined by either clause of Proposition~\ref{P.ssa}. 

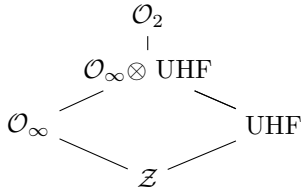
\begin{tikzpicture}
  \matrix[row sep=0.2cm,column sep=.2cm] {
& \node (O2)  {$\cO_2$}; & \\
&  \node  (OinftyUHF) {$\cO_\infty\otimes$ UHF}; 
\\
  \node  (Oinfty) {$\cO_\infty$}; &  & 
  \node (UHF) {UHF}; \\
& \node (Z) {$\cZ$}; 
\\
};
\path (OinftyUHF) edge [-]   (O2); 
\path (UHF) edge [-]   (OinftyUHF); 
\path (Oinfty) edge [-]   (OinftyUHF); 
\path (UHF) edge [-]   (OinftyUHF); 
\path (Z) edge [-] (Oinfty) ; 
\path (Z) edge [-] (UHF) ; 
\end{tikzpicture}

Finding an s.s.a. algebra other than the ones in the diagram 
would refute the revised Elliott program.

\section{Tracial von Neumann algebras} 
Many of the pathologies that plague (or enrich, depending on the point of view)
 the theory of C*-algebras are not present in von Neumann algebras. 

By a result of McDuff, the relative commutant of a II$_1$ factor in its ultrapower is trivial, 
nontrivial and abelian, or the factor  tensorially absorbs $R$ (see Proposition~\ref{P.McD}). 
 Each of these three classes is nonempty, and there is presently  no other known 
 method for distinguishing  theories of II$_1$ factors (see \cite{FaHaSh:Model3}). 

The hyperfinite II$_1$ factor $R$  is a canonical object and every embedding of $R$ into a model of its theory 
is elementary (\S\ref{S.s.s.a..model}).  
However, there are embeddings between models of the theory of $R$ that are not elementary
(i.e., the theory of $R$ is  not \emph{model-complete}), 
and in particular this theory  does not allow  the elimination of quantifiers 
 (\cite{FaGoHaSh:Existentially}, \cite{GoHaSi:Theory}).  This may be an indication that we do not  have the right language for the 
 theory of II$_1$ factors. The obstruction for the elimination of quantifiers extracted in \cite{FaGoHaSh:Existentially} 
 from \cite{Jung:Amenability} is removed by adding a predicate for the unitary conjugacy relation.  
As this is a definable relation, adding such predicate affects only syntactical structure of the language. 
 It is not clear whether adding finitely, or even countably, many such predicates could make the 
 theory of $R$ model-complete. 
This may suggest  that the theory of $R$ is as complicated as the first-order arithmetic or ZFC.

Given a II$_1$ factor $M$ and a projection $p$ in $M$, are $M$ and its corner $pMp$ elementarily equivalent? 
By the Keisler--Shelah theorem, this is equivalent to asking whether these algebras have isomorphic ultrapowers. 
A positive answer would imply that all free group factors $L(F_n)$, for $n\geq 2$,  are elementarily equivalent, 
giving a `poor man's' solution to the well-known problem whether the free group factors are isomorphic (see \cite{dykema1994interpolated}). 
On the other hand, a negative answer would provide a continuum of distinct theories of~II$_1$ factors
that are corners of  $L(F_2)$ . 
A deeper analysis of the model theory of II$_1$ factors will nceessarily involve  Voiculescu's free probability.

In recent years theories of C*-algebras and von Neumann algebras are increasingly considered as inseparable. 
Some of the most exciting progress on understanding tracial C*-algebras was initiated  in  \cite{MaSa:Strict}. 
We shall return to this in \S\ref{S.Relative}, but see also \cite{Bro:Symbiosis}.

\section{Massive algebras I: Ultraproducts} 
\label{S.Ultra} 

We  now consider algebras that are rarely nuclear
  and never  separable, but are nevertheless indispensable tools in the study of separable nuclear algebras. 

\emph{Ultraproducts} emerged independently in logic and in functional analysis (more precisely, in the theory of 
II$_1$ factors) in the 1950's (see the introduction to \cite{She:Notes}). 
If $(A_n,d_n)$, for $n\in \bbN$, are bounded metric structures of the same signature and $\cU$ is an ultrafilter 
on $\bbN$, then the \emph{ultraproduct} $\prod_{\cU} A_n$ is defined as follows. 
On the product structure $\prod_n A_n$ consider the quasi-metric 
\[
\textstyle d_{\cU}((a_n), (b_n))=\lim_{n\to \cU} d_n(a_n,b_n). 
\]
Since every function symbol $f$ has a fixed modulus of uniform continuity, it defines a 
uniformly continuous function on the quotient metric structure $\prod_n A_n/\sim_{d_{\cU}}$. 
This structure is the \emph{ultraproduct} of $A_n$, for $n\in \bbN$,
associated to the ultrafilter $\cU$. It is denoted by~$\prod_{\cU} A_n$.

In the not necessarily bounded case one replaces  $\prod_n A_n$ 
with $\{(a_n)\in \prod_n A_n: a_n$ belong to the same domain of quantification$\}$. 
With our conventions, in the operator algebra case this is 
the $\ell_\infty$-product usually denoted $\prod_n A_n$. The nontrivial fact that an ultrapower of tracial von Neumann algebras
is a tracial von Neumann algebra is an immediate consequence of the axiomatizability.

The usefulness of ultraproducts  
draws its strength largely from  two basic principles. The first one is \emph{\L o\'s's theorem}, 
stating that for any formula $\phi(\bar x)$ we have 
\[
\textstyle\phi(\bar a)^{\prod_{\cU} A_n}=\lim_{n\to \cU} \phi(\bar a_n)^{A_n}. 
\]
This in particular implies that the diagonal embedding of $A$ into its ultrapower is  elementary (\S\ref{S.ES}), 
and therefore the theory is preserved by taking ultrapowers. 
The second principle will be discussed in \S\ref{S.CtbleSat}.

This may be a good place to note two results in abstract model theory that carry over to the metric case (\cite{BYBHU}).  
A  category $\bbK$ with an  appropriately defined ultraproduct construction
is closed under the elementary equivalence
 if and only if it is closed under isomorphisms, ultraproducts, and ultraroots (i.e., $A^{\cU}\in \bbK$ implies $A\in \bbK$).  
By the Keisler--Shelah theorem, two models are elementarily equivalent if and only if they 
have isomorphic ultrapowers. Both results require considering ultrafilters on arbitrarily large sets
(see \cite{shelah1992vive}).

The fact that it is easier to prove that an ultraproduct of C*-algebras is a C*-algebras 
than that an ultraproduct of tracial von Neumann algebras is a tracial von Neumann algebra
is reflected in the fact that it is easier to prove that the C*-algebras are axiomatizable  than 
that the tracial von Neumann algebras are axiomatizable.

All ultrafilters considered here concentrate on $\bbN$ and are nonprincipal. 
It is not possible to construct such an ultrafilter in ZF alone, as a (rather weak) form of
the Axiom of Choice is required for its construction. However, results about separable C*-algebras and
separably acting II$_1$ factors proved using ultrafilters can be proved without appealing to 
the Axiom of Choice, by standard absoluteness arguments.

An ultrapower of an infinite-dimensional, simple, unital C*-algebra is by \L o\'s's theorem
unital.  It is, however, nonseparable, not nuclear, and it is simple only under  exceptional circumstances.  
This shows that none of these three properties is axiomatizable (cf. Proposition~\ref{P.4.1}).   
Nevertheless, 
\snus{} C*-algebras can be constructed by using the Henkin construction and omitting types theorem (\cite{FaMa:Omitting}, 
see \S\ref{S.Omitting}). 

\vskip.2in

\subsubsection{Countable saturation} \label{S.CtbleSat}

We define the second important  property of massive algebras.
If a type 
(see \S\ref{S.Omitting}) 
is allowed to contain formulas with parameters from 
an algebra $A$ we say that it is a type \emph{over}~$A$. 

An algebra $A$ is \emph{countably saturated} if  every countable type $\bft(\bar x)$ over $A$ 
is realized in $A$ if and only if it is consistent. 
(These algebras are sometimes said to be \emph{$\aleph_1$-saturated}. The latter terminology is  more conveniently 
extended to higher cardinalities.)
Every ultrapower associated to a nonprincipal ultrafilter on $\bbN$ is countably saturated. 
A weakening of countable saturation suffices for many purposes (see \S\ref{S.Coronas}), 
and we shall return to full saturation in \S\ref{S.Automorphisms}.

\subsection{Relative commutants} 
\label{S.Relative} 
In the theory of operator algebras 
even more important than the ultrapower itself is the \emph{relative commutant} of the algebra
inside the ultrapower, 
\[
A'\cap A^{\cU}=\{b\in A^{\cU}: ab=ba\text{ for all }a\in A\}. 
\]
The current prominence of ultrapowers as a tool for studying separable algebras can be traced back 
to  McDuff (\cite{McDuff:Central}) and the following proposition (generalized to s.s.a. algebras in \cite{ToWi:Strongly}).

\begin{proposition} \label{P.McD} 
If $D$ is strongly self-absorbing and $A$ is separable, then 
$A$ is $D$-absorbing if and only if $D$ embeds into $A'\cap A^{\cU}$. 
\end{proposition} 

The nontrivial, converse, implication  uses  the following (a lemma in model theory that I learned from
Wilhelm Winter)  proved using the intertwining argument. 

\begin{lemma} If $A\subseteq B$ are separable metric structures and $B^{\cU}$ has a sequence of isometric automorphisms $\alpha_n$
such that $\lim_n \alpha_n(a)=a$ for all $a\in A$ and $\lim_n \dist(\alpha_n(b), A)=0$ for all $b\in B$, then $A$ and $B$ are isometrically isomorphic. 
\end{lemma}

Noting that all nonprincipal ultrafilters on $\bbN$ 
`look the same' and in particular that  the choice of $\cU$ in  Proposition~\ref{P.McD}  is irrelevant as long as it is a nonprincipal ultrafilter 
on $\bbN$, one may ask the following. 

\begin{question} \label{Q.Ultrapower} 
 If   $M$ is a separable metric structure, does the isomorphism type 
 of~$M^{\cU}$ (and $M'\cap M^{\cU}$, if $M$ is a Banach algebra)  depend on $\cU$ at all? 
\end{question} 

If $M$ is a Hilbert space or a measure algebra, then a simple argument (using 
Maharam's theorem in the latter case) gives a negative answer. 
Also, Continuum Hypothesis (CH) implies negative answer to both questions for an arbitrary 
separable~$M$  (see \S\ref{S.Automorphisms}). 
Therefore, the question is whether CH can be removed from this proof. 

Question~\ref{Q.Ultrapower} for relative commutants 
was  asked by McDuff (\cite{McDuff:Central}) 
and Kirchberg (\cite{Kirc:Central}) in the case of McDuff  factors and C*-algebras, respectively. 
In  \cite{GeHa} it was  proved that, under some additional 
assumptions on $M$, CH is equivalent to the positive answer to either of these questions!
 This was achieved by using only results from classical (`discrete') model theory. 
By using the logic of metric structures and Shelah's non-structure theory, 
 the full result was proved in  \cite{FaHaSh:Model1} and \cite{FaSh:Dichotomy}. 

\begin{theorem} Assume CH fails. 
If $M$ is a separable C*-algebra or a McDuff factor with a separable 
predual, then $M$ has  $\tfc$ nonisomorphic ultrapowers and $\tfc$ nonisomorphic
relative commutants associated to nonprincipal ultrafilters on $\bbN$. 
\end{theorem}

Let's zoom out a bit. 
A complete first-order theory $\bfT$ has the \emph{order property} if there exist $n\geq 1$ and 
a $2n$-ary formula $\phi(\bar x, \bar y)$ such that for every $m$ there is a model 
$\fM$ of $\bfT$ which has a `$\phi$-chain' of length at least $m$. 
A \emph{$\phi$-chain} is a sequence $\bar x_i, \bar y_i$, for $i\leq m$, 
such that 
\[
\phi(\bar x_i,\bar y_j)=0\text{ if }i\leq j 
\text{ and } 
\phi(\bar x_i,\bar y_j)=1\text{ if }i>j. 
\]
This is the metric version  of one of the important non-structural properties of theories in Shelah's stability theory (\cite{She:Classification} 
and  \cite{FaHaSh:Model1}). 
The theory of any infinite-dimensional C*-algebra and of any II$_1$ factor has the order property. 
This is proved by continuous functional calculus and by utilizing noncommutativity, respectively.   
However, the theories of abelian tracial von Neumann algebras do not have the order property, essentially by applying Maharam's theorem
on measure algebras.

\begin{theorem}\label{T.Dichotomy}
Suppose that $A$ is a separable structure  in a separable
language.
\begin{enumerate}
\item\label{T0.1}
If the theory of $A$ does not have the order property then all of its ultrapowers associated
to nonprincipal ultrafilters on $\bbN$ are isomorphic.  
\item \label{T0.2} If the theory of $A$ has the order property then the following are equivalent:
\begin{enumerate}
\item $A$ has  fewer than $2^{2^{\aleph_0}}$ nonisomorphic ultrapowers associated with
nonprincipal ultrafilters on $\bbN$.
\item 
all  ultrapowers of $A$  associated
to nonprincipal ultrafilters on $\bbN$ are isomorphic.  
\item the Continuum Hypothesis holds.
\end{enumerate}
\pushcounter
\end{enumerate}
\end{theorem}

\subsection{Model theory of the relative commutant} 
The notion of a relative commutant does not seem to have  a useful generalization 
in the  abstract model theory and its model-theoretic properties are still poorly understood. 

While the structure of relative commutants of II$_1$ factors in their ultrapowers provides the only known method for distinguishing 
their theories, 
every infinite-dimensional separable C*-algebra has a nontrivial relative commutant in 
 its ultrapower (\cite{Kirc:Central}, also \cite{FaHaSh:Model1}). 
 The relative commutant of the Calkin algebra (\S\ref{S.Coronas}) in its ultrapower is 
 trivial (\cite{Kirc:Central}) and the relative commutant of $\cB(H)$ may or may not 
 be trivial, depending on the choice of the ultrafilter (\cite{FaPhiSte:Relative}).

It is not difficult to see that the existential theory of $A'\cap A^{\cU}$ depends only on the theory of $A$. 
However, a result of \cite{Kirc:Central} implies that there is a separable C*-algebra $A$ elementarily equivalent to $\cO_2$ 
such that $A'\cap A^{\cU}$ and $\cO_2\cap \cO_2^{\cU}$ have different \aet-theories. 
(An \emph{\aet-sentence} is one of the form $\sup_{\bar x}\inf_{\bar y} \phi(\bar x, \bar y)$ where $\phi$ is quantifier-free.) 
In the following all ultrafilters are nonprincipal ultrafilters on~$\bbN$. 

\begin{proposition} \label{P.rc} Assume $A$ is a separable C*-algebra. 
\begin{enumerate}
\item \label{I.P.rc.1} For all $\cU$ and $\cV$, the algebras $A'\cap A^{\cU}$ and $A'\cap A^{\cV}$ are elementarily equivalent. 
\item \label{I.P.rc.2} For every separable $C\subseteq A'\cap A^{\cU}$ we have $\ThE(A'\cap C'\cap A^{\cU})=\ThE(A'\cap A^{\cU})$. 
\item \label{I.P.rc.3} If $D$ is a separable unital subalgebra of $A'\cap A^{\cU}$ then there are $\aleph_1$ commuting 
copies of $D$ inside $A'\cap A^{\cU}$. 
\end{enumerate}
\end{proposition}

An entertaining  proof of \eqref{I.P.rc.1} can be given by using basic set theory. 
Collapse $2^{\aleph_0}$ to $\aleph_1$ without adding reals. 
Then $\cU$ and $\cV$ are still  ultrafilters on $\bbN$ and  one can use saturation to find an isomorphism 
between the ultrapowers that sends $A$ to itself. The theories of two algebras are unchanged, and therefore by absoluteness
the result follows. 
Clause \eqref{I.P.rc.3} is an immediate consequence of \eqref{I.P.rc.2} and it is a minor strengthening of a result in \cite{Kirc:Central}.  

When $A$ is not $\cZ$-stable, the relative commutant of $A$  can  have characters even if it is simple (\cite{KirRo:Central}).  
In the case when algebra $A$ is nuclear and $\cZ$-stable, $A'\cap A^{\cU}$ inherits some properties from $A$. 
For example, each of the  traces on $A'\cap A^{\cU}$  extends to a trace on $A^{\cU}$ by \cite{MaSa:Strict} (cf. Proposition~\ref{P.Elementary}). 
The relative commutants of s.s.a. algebras 
are well-understood; the following was proved in \cite{FaHaRoTi}.

\begin{proposition} If $D$ is a s.s.a. algebra and $\cU$ is a nonprincipal ultrafilter on $\bbN$, 
then $D'\cap D^{\cU}$ is an elementary submodel of $D^{\cU}$. 
Moreover, CH implies that these two algebras are isomorphic. 
\end{proposition} 

\subsection{Expansions and traces} If a metric structure $A$ is expanded by adding a new  
predicate $\tau$, its ultrapower $A^{\cU}$  
expands to the ultrapower of the expanded structure $(A,\tau)^{\cU}$ which still satisfies \L o\'s's theorem and is countably saturated. 

If $A$ is a unital tracial C*-algebra then its traces form a weak*-compact convex subset $T(A)$ of the dual unit ball. 
For $\tau\in T(A)$ denote the tracial von Neumann algebra associated with the  $\tau$-GNS representation  (\S\ref{S.OA}) by $N_\tau$. 
If $A$ is simple and infinite-dimensional and $\tau$ is an extremal trace  then $N_\tau$ is a factor, and if $A$ is in addition nuclear and separable 
then $N_\tau$ is isomorphic to 
the hyperfinite factor $R$. This is because $A$ is nuclear if and only if its weak closure in 
every representation is an injective von Neumann algebra, and $R$ is the only injective II$_1$ factor with a separable predual.  
The following was proved in \cite{MaSa:Strict} and improved to the present form in \cite{KirRo:Central}. 

\begin{proposition} \label{P.MS} If $A$ is separable and $\tau\in T(A)$, then 
the quotient map from $A'\cap A^{\cU}$ to $N_\tau'\cap (N_\tau)^{\cU}$ is surjective. 
\end{proposition} 

If $b\in A^{\cU}$ is such that its image is in the commutant of $N_\tau'$, 
then by countable saturation one finds a positive element $c$ of norm 1 such that $\tau(c)=0$ and 
$c(a_nb-ba_n)=(a_n b - b a_n)c=0$ for all $a_n$ in a fixed countable dense subset of $A$. The fact that the type of such $c$ is
consistent follows from the fact that the image of $b$ is in $N_\tau'$. 
Then $(1-c)b(1-c)$ is in $A'\cap A^{\cU}$ and it has the same image under the quotient map as $b$. 

Proposition~\ref{P.MS}  precipitated remarkable progress on understanding tracial C*-algebras, 
the most recent results of which are \cite{matui2013decomposition} and \cite{sato2014nuclear}.

\section{Massive algebras II: Coronas} \label{S.Coronas} 
Another class of massive C*-algebras (with no analogue in von Neumann algebras) 
has special relevance to the study of separable algebras. 
If $A$ is a non-unital C*-algebra, the \emph{multiplier algebra} of $A$, $M(A)$, is the 
noncommutative analogue of the \v Cech--Stone compactification of a locally compact Hausdroff space. 
It is the surjectively universal unital algebra containing $A$ as an essential ideal. 
The \emph{corona} (or \emph{outer multiplier}) algebra of $A$ is the quotient  $M(A)/A$. 
Some examples of coronas are the Calkin algebra $\cQ(H)$ (the corona of the algebra of compact operators) and 
the \emph{asymptotic sequence algebra} $\ell_\infty(A)/c_0(A)$ for a unital $A$. 
The latter algebra, as well as the associated \emph{central sequence algebra} $A'\cap \ell_\infty(A)/c_0(A)$
are sometimes used in classification of C*-algebras 
instead of the metamathematically heavier ultrapowers and the corresponding relative
commutants. 
While  \L o\'s's theorem miserably fails for the asymptotic sequence algebra, all coronas and corresponding 
relative commutants share some properties of  countably saturated algebras. 
The simplest of these properties is being SAW*: 
for any two  orthogonal separable subalgebras $A$ and $B$ of a corona there exists a positive element $c$
such that $ca=a$ for all $a\in A$ and $cb=0$ for all $b\in B$. 

\vskip .2in

\subsubsection{Quantifier-free saturation} An algebra $C$ is \emph{quantifier-free saturated} if every countable type over  $C$ 
consisting only of quantifier-free formulas  is consistent if and only if it is realized in $C$. 
An algebra $C$ is \emph{countably degree-1 saturated} if every countable type over  $C$ 
consisting only of formulas of the form $\|p\|$, 
where $p$ is a *-polynomial of degree 1,  is consistent if and only if it is realized in $C$. 
A dummy variable argument shows that the degree-2 saturation is equivalent to quantifier-free saturation. 
By refining an argument introduced by Higson, the following was proved in \cite{FaHa:Countable}. 

\begin{theorem} If $A$ is a corona of a separable non-unital C*-algebra, or a relative commutant
of a separable subalgebra of such corona, then $A$ is countably degree-1 saturated. 
\end{theorem}

A very interesting class of countable degree-1 saturated C*-algebras was isolated in \cite{Voi:Countable}. 

\subsubsection{A sampler of properties of countable degree-1 saturated algebras} 
Assume $C$ is  
countably degree-1 saturated (the results below also apply to tracial von Neumann algebras, and in this case 
 (1), (3) and (5) do not even require
countable degree-1 saturation). 

(1) $C$ has SAW* as well as every other known countable separation property~(\cite{FaHa:Countable}).

(2) A separable algebra $A$ is isomorphic to a unital subalgebra of $C$ if and only if $\ThE(A)\subseteq\ThE(C)$.

(3) A representation of a group $\Gamma$ in $A$  is a homomorphism 
$\pi\colon \Gamma\to (\GL(A),\cdot)$. It is
 \emph{unitarizable} if there is  an invertible $h\in A$ such that $h^{-1}\pi(g)h$ is 
a unitary for all $g\in \Gamma$.
Conjecturally unitarizability of all uniformly bounded representations of a group $\Gamma$ on $\cB(H)$  is equivalent to  
the  amenability of $\Gamma$ (see \cite{pisier2001similarity}). 
If $\Gamma$ is a countable amenable group, then every uniformly bounded representation 
$\pi$ of $\Gamma$ in $C$ is unitarizable (\cite{choi2013nonseparable}). 

(4) $C$ is not isomorphic to the tensor product of two infinite-dimensional 
  algebras (\cite{fang2006central} for the ultraproducts of II$_1$ factors 
  and \cite{Gha:SAW*} for the general result). Therefore  
an ultrapower or a corona is never 
isomorphic  to a nontrivial tensor product and   the separability assumption is needed in Proposition~\ref{P.McD}.

(5) (`Discontinuous  functional calculus.')
If $a$ is a normal operator, 
then by the \emph{continuous  functional calculus} for every continuous complex-valued function $g$ on 
the spectrum, 
$\Sp(a)$, of $a$ the naturally defined $g(a)$ belongs to the abelian algebra generated by $a$. 
 
\begin{proposition} \label{P.CD1} 
 Assume $C$ is countably degree-1 saturated and 
 $B\subseteq \{a\}'\cap C$ is separable, $U\subseteq \Sp(a)$ is open, 
and $g\colon U\to \bbC$ is a bounded continuous function. 
Then there exists $c\in C\cap C^*(B, a)'$ such that for every  
$f\in C_0(\Sp(a))$   we have 
\[
cf(a)=(gf)(a). 
\]
If moreover $g$ is real-valued then $c$ can be chosen to be self-adjoint. 
  \end{proposition} 

The `Second Splitting Lemma' (\cite[Lemma~7.3]{BrDoFi:Unitary})  is a special case of the above  when $C$ is 
the Calkin algebra, $a=h_0$ is self-adjoint,   
 and the range of $g$ is $\{0,1\}$. 
 
\vskip .2in
  
\subsubsection{Failure of saturation}\label{S.Failure} 
  While the asymptotic sequence algebras, as well as some abelian coronas, are fully countably saturated (\cite{FaSh:Rigidity}), 
this is not true for sufficiently noncommutative coronas. 
By a K-theoretic argument N. C. Phillips constructed two unital embeddings of the CAR algebra into the Calkin 
algebra $\cQ(H)$ that are approximately unitarily equivalent, but not conjugate by a unitary (\cite[\S 4]{FaHa:Countable}). 
This gives a countable quantifier-free type over $\cQ(H)$  that is consistent but not realized. 
Even coronas of separable abelian C*-algebras provide a  range of different saturation 
properties  (see \cite{FaSh:Rigidity}). 

\vskip.2in
\subsection{Automorphisms} \label{S.Automorphisms} 
A metric model $A$ is \emph{saturated} if every type over $A$ whose cardinality is smaller than the \emph{density character} $\chi(A)$ of $A$
(i.e., the smallest cardinality of a dense subset) 
 which is consistent is realized in $A$.  The Continuum Hypothesis (CH) 
implies that  all countably saturated models of cardinality $2^{\aleph_0}$ are saturated. 
A transfinite back-and-forth argument shows 
  that any two elementarily equivalent saturated models of the same density character are isomorphic and 
 that a saturated model $A$ has $2^{\chi(A)}$ automorphisms. By a counting argument, most of these automorphisms 
 are outer and moreover nontrivial when `trivial automorphism' is defined in any reasonable way; 
 see \cite{CoFa:Automorphisms} for a (lengthy) discussion. 
 This explains the effectiveness 
of CH  as a tool for resolving problems of a certain form. 
A deeper explanation is given in  Woodin's celebrated $\Sigma^2_1$-absoluteness theorem (see \cite{Wo:Beyond}). 

By the above, CH implies that an ultrapower $A^{\cU}$ of a separable, infinite-dimensional algebra  has automorphisms
that do not lift to automorphisms of $\ell_\infty(A)$. 
Much deeper is a complementary series of results of Shelah, to the effect that if  ZFC is consistent then 
so is the assertion that any isomorphism between ultraproducts of models with the strong independence property 
lifts to an isomorphism of the products of these models~(\cite{shelah2008vive}). No continuous version of this
result is known. 
One difficulty in taming ultrapowers is that the ultrafilter is not a definable object; in particular 
Shelah's results apply only to a carefully constructed ultrafilter in a specific model of ZFC. 

Motivated by work on extension theory and a very concrete question about the unilateral 
shift, in \cite{BrDoFi:Unitary} it was asked whether the Calkin algebra has outer automorphisms. 
Since the Calkin algebra is not countably saturated (\S\ref{S.Failure}) 
it took some time before such an automorphism was constructed 
using CH   
(\cite{PhWe:Calkin}). This is one of the most complicated known CH constructions, 
involving an intricate use of EE-theory to extend isomorphisms of direct limits of separable subalgebras. 
A simpler proof was given in \cite[\S 1]{Fa:All}, and the method was further refined in \cite{CoFa:Automorphisms}. 
Instead of following the usual back-and-forth construction in which isomorphisms between separable 
subalgebras are recursively extended, one uses CH to embed the first derived limit of an inverse system of abelian groups 
into the outer automorphism group. 

Forcing axioms imply that the Calkin algebra has only inner automorphisms (\cite{Fa:All}). 
Conjecturally, for every non-unital separable C*-algebra the assertion that its corona has only (appropriately defined) 
`trivial' automorphisms is independent of ZFC (see \cite{CoFa:Automorphisms}). Even the abelian case of this conjecture
is wide open (\cite{FaSh:Rigidity}). 

The `very concrete question' of 
 Brown--Douglas--Fillmore alluded to two paragraphs ago  is still wide open: Is there an automorphism of $\cQ(H)$ 
 that sends the image of the unilateral shift $\dot s$  to its adjoint? 
 Fredholm index obstruction shows that such an automorphism cannot be inner. Since the nonexistence
 of outer automorphisms of $\cQ(H)$ is relatively consistent with ZFC,  so is a  negative 
 answer to the BDF question. Every known  automorphism $\alpha$ of $\cQ(H)$ in every model of ZFC
 has the property that its restriction to any separable subalgebra is implemented by a unitary. 
 Both $\dot s$ and $\dot s^*$ are unitaries with full spectrum and  
 no nontrivial roots. 
 It is, however,  not even known whether $\dot s$ and $\dot s^*$  have the same (parameter-free) type in $\cQ(H)$; 
 a positive answer would provide a strong motivation for the question of whether $\cQ(H)$ is countably homogeneous.

\vskip.2in

\subsection{Gaps} \label{S.Gaps} 
A \emph{gap} in a semilattice $\cB$  is a pair $\cA$, $\cB$ such that 
$a\wedge b=0$ for all $a\in \cA$ and all $b\in \cB$ but there is 
no $c$ such that $c\wedge a=a$ and $c\wedge b=0$ for all $a\in \cA$ and $b\in \cB$. 
There are no countable gaps in a countably saturated Boolean algebra such as $\PNF$, the quotient of $\cP(\bbN)$ over 
the ideal $\Fin$ of finite sets. 
In 1908 Hausdorff constructed a gap in $\PNF$ with both of its sides of cardinality $\aleph_1$.  
Later Luzin constructed 
a family of $\aleph_1$ orthogonal elements in $\PNF$ such that any two of its disjoint uncountable 
subsets form a gap.  
It should be emphasized that both results were proved 
without using CH or any other additional set-theoretic axioms. 

Hausdorff's and Luzin's results  show that $\PNF$ is not more than countably saturated. 
In particular, if the Continuum Hypothesis fails then the obvious back-and-forth method for constructing automorphisms
of $\PNF$ runs into difficulties after the first $\aleph_1$ stages. 
In one form or another, gaps were used as an obstruction to the existence of morphisms in several 
consistency results in analysis, notably as obstructions to extending 
a partial isomorphism (\cite[\S V]{Sh:Proper}, \cite{DaWo:Introduction}, \cite{Fa:All}).

Two subalgebras $A$ and $B$ of an ambient algebra $C$ form 
a \emph{gap} if $ab=0$ for all $a\in A$ and $b\in B$, but there is no 
positive element $c$ such that $ca=a$ and $cb=0$ for all $a\in A$ and all $b\in B$. 
The gap structure of $\PNF$ can be imported into the Calkin algebra, but
the gap structure of the latter is also much richer~(\cite{Za-Av:Gaps}). 

However, the failure of higher saturation in coronas is also manifested in a genuinely noncommutative fashion. 
A countable family of commuting operators in a corona of a separable algebra can be lifted to a 
family of commuting operators if and only if this is true for each one of its finite subsets.

\begin{proposition} \label{P.Twist} 
In $M_2(\ell_\infty/c_0)$ there exists a family 
of $\aleph_1$ 
orthogonal projections such that none of its uncountable subsets can be lifted to a commuting family of projections in $M_2(\ell_\infty)$. 
\end{proposition} 

This was stated in \cite{FaWo:Set} for the Calkin algebra in place of (barely noncommutative) $M_2(\ell_\infty/c_0)$, but the proof given there clearly gives the stronger result. 
The combinatorial essence for the proof of  Proposition~\ref{P.Twist}  echoes Luzin's original idea. One recursively 
constructs projections $p_\gamma$
in $M_2(\ell_\infty)$ so that  $p_\gamma p_{\gamma'}$ is compact 
but $\|[p_\gamma,p_{\gamma'}]\|>1/4$ for all $\gamma\neq\gamma'$. 
Then the image this family in the corona is as required, as 
a counting argument  shows that no uncountable 
subfamily can be simultaneously diagonalized.

Recall that every uniformly bounded representation of a countable amenable group in a
countably degree-1 saturated algebra is unitarizable (Proposition~\ref{P.CD1}).
This is false for uncountable groups. This was proved in \cite{choi2013nonseparable}  and improved to the present form  in
  \cite{vignati2014algebra} using Luzin's gap.

\begin{proposition} \label{P.Gamma}There is a uniformly bounded representation $\pi$ of 
$\bigoplus_{\aleph_1}\bbZ/2\bbZ$ 
on $M_2(\ell_\infty/c_0)$ such that the restriction of $\pi$ to a subgroup is unitarizable if and only if the subgroup is 
countable. 
\end{proposition} 

The construction of Kadison--Kastler-near, but not isomorphic, 
nonseparable algebras in \cite{choi1983completely}  
involves what at the hindsight can be considered as a gap. 
It is not known whether there is a separable example  (see \cite{CSSWW:Perturbations}
for several partial positive results). 

\vskip.1in
\section{Nonseparable algebras} 
Not surprisingly, the theory of nonseparable algebras hides surprises and problems not present in the separable 
case; see \cite{We:Prime}.

\vskip.2in
\subsection{Nonseparable UHF algebras} 
Uniformly hyperfinite (UHF) algebras are defined as tensor products of full matrix algebras (\S\ref{S.s.s.a.}). 
However, there are two other natural ways to define uniformly hyperfinite: as (i) an inductive limit of a net 
of full matrix algebras, or (ii) as an algebra in which every finite subset can be arbitrarily well approximated 
by a full matrix subalgebra. These three notions, given in the order of decreasing strength, 
 coincide in the separable unital case.  Dixmier asked whether 
separability is needed for this conclusion.    The answer is that in every uncountable density character,  
UHF and (i) differ, but that one needs an algebra of density character  $\aleph_2$ in order to distinguish 
between (i) and (ii) (\cite{FaKa:Nonseparable}). 
An extension of methods of \cite{FaKa:Nonseparable} 
 resulted in a nuclear, simple C*-algebra that has irreducible representations on both separable and 
nonseparable Hilbert space (\cite{Fa:Graphs}). This is in contrast with the transitivity of the 
space of irreducible representations of a separable simple C*-algebra (\cite{KiOzSa}).

\subsection{Representation theory}
Representation theory of separable algebras has deeply affected  development of the classical descriptive set theory, 
as evident from the terminology of both subjects (terms `smooth' and `analytic' have the same, albeit nonstandard in other 
areas of mathematics, meaning). 
Extension of  the work of Glimm and Effros  on representation theory combined with methods from logic initiated the 
abstract classification theory (\S\ref{S.Abstract}). 
The representation theory of nonseparable algebras was largely abandoned because  some of the central problems
proved to be intractable (see the introduction to \cite{AkeWe:Consistency}). 
One of these stumbling blocks,  \emph{Naimark's problem},
 was partially solved  in \cite{AkeWe:Consistency} (see also \cite{We:Set}). 
By using a strengthening of CH (Jensen's~$\diamondsuit_{\aleph_1}$ principle) and a deep result
on representation theory of separable C*-algebras (an extension of \cite{KiOzSa} mentioned above),  
Akemann and Weaver constructed
a C*-algebra that has a unique (up to spatial equivalence) irreducible representation on a Hilbert space, 
but is not isomorphic to the algebra of compact operators on any Hilbert space. 
An extension of \cite{AkeWe:Consistency} shows that $\diamondsuit_{\aleph_1}$ implies the existence of 
a simple C*-algebra with exactly $m$ inequivalent irreducible representations. By a classical result of Glimm
(closely related to the Glimm--Effros dichotomy), a simple separable C*-algebra with two inequivalent representations has  
$2^{\aleph_0}$ inequivalent  representations. 
It is not known whether a counterexample to Naimark's problem 
can be found in ZFC alone or by using  an   axiom other than $\diamondsuit_{\aleph_1}$ (such as 
$\diamondsuit_\kappa$ for $\kappa>\aleph_1$). The fact that every forcing notion that adds a new real number 
destroys all ground-model examples is a bit of an annoying teaser. 

Cyclic representations  of C*-algebras are, via the GNS construction (\S\ref{S.GNS}),  in a natural bijective correspondence with their states (i.e., 
positive unital functionals). Pure (i.e., extremal) states are noncommutative versions of ultrafilters. The space of nonprincipal ultrafilters on $\bbN$, 
 (along  with the associated quotient structure $\cP(\bbN)/\Fin$)  is arguably the most 
important  set-theoretically malleable object known to man. The study of pure states on $\cB(H)$ (i.e., `quantized ultrafilters') has already 
produced some surprising results (\cite{AkeWe:B(H)}, \cite{bice2011filters}; also see \cite{marcus2013interlacing}).

\subsection{Amenable operator algebras} 
A prominent open problem in the theory of operator algebras is whether every 
algebra of operators on a Hilbert space which is amenable is isomorphic to a C*-algebra. 
By using Proposition~\ref{P.Gamma}, one obtains the 
following (\cite{choi2013nonseparable}, \cite{vignati2014algebra}).

\begin{theorem} \label{T.amenable} There exists a nonseparable amenable  subalgebra of $M_2(\ell_\infty)$    
which is not isomorphic to a C*-algebra. None of its nonseparable 
amenable subalgebras is isomorphic to a C*-algebras, yet it is an 
inductive limit of separable subalgebras (even elementary submodels) each of which is isomorphic to a C*-algebra. 
Moreover, for every $\eps>0$ such an algebra can be found in an $\eps$-Kadison--Kastler neighbourhood of a C*-algebra.\end{theorem}

The question whether there exists a separable counterexample remains open; see \cite{marcoux2013abelian}.

\section{Concluding remarks} 
 
 The most recent wave of applications of logic to operator algebras started  
 by work of Nik Weaver and his coauthors, in which several long-standing problems were solved 
 by using additional set-theoretic axioms (see \cite{We:Set}). 
 Although we now know that the answers to some of those problems  (such as the existence
 of outer automorphisms of the Calkin algebra) are independent from ZFC, 
statements of  many prominent open problems in operator algebras
 are absolute between models of ZFC and therefore unlikely to be independent (see the appendix to 
 \cite{Fa:Absoluteness} for a discussion). 

Nevertheless, operator algebras do mix very well with logic. 
Jon Barwise said ``As logicians, we do our subject a disservice by convincing others that the logic is first-order and then convincing them that almost none of the concepts of modern mathematics can really be captured in first-order logic.''
Remarkably, some of the deepest results on the structure of C*-algebras have equivalent formulation in 
 the language of (metric)  first-order logic  (this applies e.g., to \cite{Win:Decomposition} and  
 \cite{Win:Nuclear}). 

In many  of the developments presented here methods from logic were blended with highly nontrivial  operator-algebraic methods. 
Good examples are the proof that the theory of $R$ does not allow elimination of quantifiers (\cite{GoHaSi:Theory}) 
the key component of which comes from \cite{Bro:Topological},   the already mentioned use of \cite{Jung:Amenability}, 
and blending of $\diamondsuit_{\aleph_1}$ with the transitivity of pure state space of separable
simple algebras (\cite{KiOzSa}) in \cite{AkeWe:Consistency}. 

Finally, some results in pure logic were motivated by work on operator algebras. 
Examples are Theorem~\ref{T.Dichotomy}, which is new even for discrete structures, and 
negative and positive results on omitting types  (\S\ref{S.Omitting}). 

\providecommand{\bysame}{\leavevmode\hbox to3em{\hrulefill}\thinspace}
\providecommand{\MR}{\relax\ifhmode\unskip\space\fi MR }
\providecommand{\MRhref}[2]{%
  \href{http://www.ams.org/mathscinet-getitem?mr=#1}{#2}
}
\providecommand{\href}[2]{#2}


\end{document}